\documentclass{article}
\usepackage[top = 2cm]{geometry}
\usepackage{amssymb}
\usepackage{amsmath}
\usepackage{array}
\usepackage{float}
\usepackage[caption = false ]{subfig}
\usepackage[final]{graphicx}
\usepackage{amsthm}

\newcommand\myeq{\mathrel{\overset{\makebox[0pt]{\mbox{\normalfont\tiny\sffamily def}}}{=}}}

\usepackage{hyperref}

\title{Jagged-time-step technique improving convergence order of Fernandez's Explicit Robin-Neumann scheme for the coupling of incompressible fluid with thin-walled structure} 

\author{Yiyi HUANG (yiyi\_huang\_me@outlook.com)}

\begin{document}
\maketitle

\begin{abstract}
  Inspired by Rybak's multiple-time-step technique, jagged-time-step technique is proposed and applied to Fernandez's \textbf{Explicit Robin-Neumann scheme}. For some instances, numerical experiments demonstrate higher convergence orders and accuracy with lower computation cost as time and space get refined.  \\

  *Notes: the work described in this article was done more than two years ago. This article is being written and extended. More numerical results, including but not limited to those at refinement $ \textit{rate} = 5 $, are to present.  The ideas described might be applicable to other algorithms. On the other hand, for easier implementation, the errors were computed in an approximate way at that time, which is to correct.
\end{abstract}


\section{Introduction}
For the coupling of incompressible fluid with thin-walled structure, \cite{rn} proposes  \textbf{Explicit Robin-Neumann scheme}. The scheme with first-order extrapolation yields unconditional stability and optimal accuracy of first-order in time.


Rybak \cite{multirate} developes a multiple-time-step technique and applies it to a decoupled scheme for coupled free flow and porous medium systems. The whole time interval is partitioned into some fine grids. Multiple fine grids constitute a coarse grid. Within each coarse time grid, the free flow solutions are computed at each fine time steps using information from the porous medium at the beginning of current coarse time grid. When it reaches the end of current coarse time grid,  the porous medium solutions are computed using information from the free flow solutions. The technique improves efficiency of computation, preserves orders of convergence and proves to be stable. Compared with algorithms that do not adopt the technique (namely the monolithic approach and decoupled scheme with single time step), the only disadvantage is it is a bit less accurate.

It is thus meaningful to design a technique that improves both efficiency and accuracy. This work is devoted to investigation of such a technique and its application to decoupled algorithms, such as Fernandez's \textbf{Explicit Robin-Neumann scheme} for the coupling of incompressible fluid with thin-walled structure. For convenience, the technique is named jagged-time-step technique and described as follows. 

The whole time interval is partitioned into some fine grids. A coarse grid consists of $ 10 $ fine grids. Both the fluid and structure are computed for multiple steps (the number of such steps are not necessarily equal to $ 10 $ ) within each coarse grid using latest information from each other. Let $ N_f $ denote the number of steps that the fluid is solved within each coarse grid, while $ N_s $ for the structure. 

The original decoupled algorithm, Fernandez's \textbf{Explicit Robin-Neumann scheme}, runs at all fine grids. Equivalently, both the fluid and structure are computed $ 10 $ steps during each coarse grid. Therefore, taking $ N_f = 10, N_s = 10 $ and applying the jagged-time-step technique to the algorithm does not make any change.

Because the fluid domain is of one dimension higher than that of the structure (in this work, for the simplified problem considered, the fluid domain is two-dimensional, while the structure is one-dimensional), it takes much more cost to compute a step of fluid than that of structure. Hence, for sake of efficiency, there should be constraints $ N_f < 10 $, and $ N_f + N_s \leq 20 $.

In what follows, the simplified problem and Fernandez's \textbf{Explicit Robin-Neumann scheme} from \cite{rn} are cited.  The jagged-time-step technique is applied to the scheme. Numerial experiments, conclusions and possible extensions follows. 

\section{The simplified problem and Fernandez's \textbf{Explicit Robin-Neumann}  scheme}
\subsection{The simplified problem}
Consider the simplified problem studied in \cite{rn} where the fluid is governed by the Stokes equations in a $d- $dimensional ($d = 2,3 $) domain $\Omega $ and the structure is assumed to be a linear thin-solid defined on a $ (d-1)- $ manifold $\Sigma $, with $ \partial \Omega = \Sigma \cup \Gamma^d \cup \Gamma^n  $. The coupled simplified problem reads: find the fluid velocity $ \textbf{u}: \Omega^f \times \mathbb{R}^+ \rightarrow \mathbb{R}^\textit{d} $, the fluid pressure $ p: \Omega^f \times \mathbb{R}^+ \rightarrow \mathbb{R} $, thestructuredisplacement $ \textbf{d} : \Sigma \times \mathbb{R}^+  \rightarrow \mathbb{R}^\textit{d} $ such that

\begin{eqnarray}  \label{eq:fluid}
\left\{
\begin{aligned}
\rho^f \partial_t \boldsymbol{u} - \textbf{div } \boldsymbol{\sigma}  ( \boldsymbol{u} ,p ) & = \textbf{0}   & \text{in \quad} & \Omega,   \\
\textbf{div} \boldsymbol{u} &= 0   & \text{in \quad}       & \Omega ,  \\
\boldsymbol{u} &= \textbf{0}  & \text{on \quad}  & \Gamma^d,   \\
\boldsymbol{\sigma} ( \boldsymbol{u} , p ) \boldsymbol{n}  &= \boldsymbol{f}^\Gamma  & \text{on \quad}  & \Gamma^n,     
\end{aligned}
\right.
\end{eqnarray}


\begin{eqnarray}  \label{eq:solid}
\left\{
\begin{aligned}
\boldsymbol{u} &= \dot{ \boldsymbol{d} }  & \text{on \quad} &\Sigma, \\
\rho^s \epsilon \partial_t  \boldsymbol{ \dot{d} } + \boldsymbol{L^ed} + \boldsymbol{L^v \dot{d} } 
  &= -\boldsymbol{\sigma (u, }  p \boldsymbol{ )n }   & \text{on \quad} &\Sigma, \\
\boldsymbol{ \dot{d} }  &=  \partial_t \boldsymbol{d}  & \text{on \quad} &\Sigma, \\
\boldsymbol{d}  &= \boldsymbol{0}  & \text{on \quad} &\partial \Sigma, 
\end{aligned}
\right.
\end{eqnarray}

with initial conditions

\[ \boldsymbol{u}(0) = \boldsymbol{u}^0, \quad \boldsymbol{d} (0) = \boldsymbol{d}^0, \quad
  \boldsymbol{ \dot{d} } (0)  = \boldsymbol{ \dot{d} }^0,  \]
where $ \rho^f $ denotes the fluid density,  $\rho^s $ thestructuredensity, $\epsilon $ thestructurethickness, $ \boldsymbol{ \dot{d} } $ thestructurevelocity , $ \boldsymbol{n} $ the exterior unit normal vector to $ \partial \Omega $, $ \boldsymbol{f^\Gamma } $ a given surface force on $ \boldsymbol{\Gamma^n } $, and 
\[ \boldsymbol{ \sigma } ( \boldsymbol{u} , p )  \myeq -p \boldsymbol{I} + 2 \mu \boldsymbol{ \varepsilon(u) } ,\quad  \boldsymbol{ \varepsilon(u) }  \myeq  \frac{1}{2} \boldsymbol{ ( \nabla u + \nabla u^T ) }, 
\] 
where $ \mu $ denotes the fluid dynamic viscosity. $ \boldsymbol{L}^e \text{and} \boldsymbol{L}^v $ stand for the elastic and viscous contributions respectively. 

\subsection{Notations}
For all the algorithms mentioned in this work, $ \tau $ denotes time step, while $ h $ stands for space discretization parameter.  

Given arbitrary variable $x $, the notation 
\begin{eqnarray} 
x^{n, \star} \myeq 
\left\{
\begin{aligned}
& 0 & \text{if \quad} & r  = 0,  \\
& x^{n-1} & \text{if \quad} & r  = 1,  \\
& 2x^{n-1} - x^{n-2} & \text{if \quad} & r = 2
\end{aligned}
\right.
\end{eqnarray}
 is used for interface extrapolations of order $ r $.

\subsection{Fernandez's Explicit Robin-Neumann  scheme}
Fernandez's \textbf{Explicit Robin-Neumann scheme} \cite{rn} is cited here. 

\par\noindent\rule{\textwidth}{1pt}
\textbf{(Fernandez) Explicit Robin-Neumann  scheme (time semi-discrete)} 
\par\noindent\rule{\textwidth}{.5pt}

1. Fluid step: find $ \textbf{u}^n: \Omega^f \times \mathbb{R}^+ \rightarrow \mathbb{R}^\textit{d} $, $ p^n: \Omega^f \times \mathbb{R}^+ \rightarrow \mathbb{R} $ such that

\begin{eqnarray}   \label{eq:ERN1fluid}
\left\{
\begin{aligned}
\rho^f \partial_\tau \boldsymbol{u}^n - \textbf{div } \boldsymbol{\sigma}  ( \boldsymbol{u}^n ,p^n ) & = \textbf{0}   & \text{in \quad} & \Omega,   \\
\textbf{div} \boldsymbol{u}^n &= 0   & \text{in \quad}       & \Omega ,  \\
\boldsymbol{u}^n &= \textbf{0}  & \text{on \quad}  & \Gamma^d,   \\
\boldsymbol{\sigma} ( \boldsymbol{u}^n , p^n ) \boldsymbol{n}  &= \boldsymbol{f}^\Gamma  & \text{on \quad}  & \Gamma^n,     \\
\boldsymbol{\sigma} ( \boldsymbol{u}^n , p^n ) \boldsymbol{ n }  + \frac{\rho^s \epsilon}{\tau} \boldsymbol{u}^n  &= \frac{\rho^s \epsilon}{\tau} \boldsymbol{\dot{d}^{n-1} } - \boldsymbol{L^ed^*} - \boldsymbol{L^v \dot{d}^* }  & \text{on \quad}  & \Sigma. \\
\end{aligned}
\right.
\end{eqnarray}

2. Solid step: find $ \textbf{d}^n : \Sigma \times \mathbb{R}^+  \rightarrow \mathbb{R}^\textit{d} $ such that

\begin{eqnarray}    \label{eq:ERN1solid}
\left\{
\begin{aligned}
\rho^s \epsilon \partial_\tau  \boldsymbol{ \dot{d}^n } + \boldsymbol{L^ed} + \boldsymbol{L^v \dot{d}^n } 
  &= -\boldsymbol{\sigma} ( \boldsymbol{u^n },   p^n ) \boldsymbol{ n }   & \text{on \quad} &\Sigma, \\
\boldsymbol{ \dot{d}^n }  &=  \partial_\tau \boldsymbol{d}^n  & \text{on \quad} &\Sigma, \\
\boldsymbol{d}^n  &= \boldsymbol{0}  & \text{on \quad} &\partial  \Sigma , 
\end{aligned}
\right.
\end{eqnarray}

\par\noindent\rule{\textwidth}{.5pt}

\vspace{2em}
Fernandez's \textbf{Explicit Robin-Neumann  scheme (time semi-discrete)} says that in each time step, first solve the fluid with Robin condition  $ (\ref{eq:ERN1fluid})_5 $ on interface with data of structure from last time step ( e.g. $ \boldsymbol{ \dot{d}^{n-1} } $ ) or by certain extrapolation strategies ( e.g. $ \boldsymbol{d^*}, \boldsymbol{\dot{d}^* } $, see section 3.1 in \cite{rn} for details ) and then solve thestructurewith Neumann condition $(\ref{eq:ERN1solid})_2 $ on interface with latest data computed from fluid. It is the Robin-Neumann conditions on interface that guarantee the stability (free of added-mass effect). With finite element discretization in space, involving variational residuals of fluid stresses on the interface,the fully discrete version of the  preceding algorithm is detailed as follows ( see \textbf{Algorithm 5} in \cite{rn} ).

\par\noindent\rule{\textwidth}{1pt}
\textbf{Fernandez Explicit Robin-Neumann  scheme (intrinsic formulation) }
\par\noindent\rule{\textwidth}{.5pt}

1. Fluid step: Find $ ( \boldsymbol{u}_h^n , p_h^n ) \in \boldsymbol{V}^f \times \boldsymbol{Q}_h $ such that

\begin{eqnarray}    \label{eq:ERN5fluid}
\left\{
\begin{aligned}
\rho^f ( \partial_\tau \boldsymbol{u}_h^n , \boldsymbol{v}_h ) + a( \boldsymbol{u}_h^n , \boldsymbol{v}_h ) +b(p_h^n, \boldsymbol{v}_h ) - b(q_h, \boldsymbol{u}_h^n ) + s_h ( p_h, q_h )   \\
+ \frac{\rho^s \epsilon}{\tau} ( \boldsymbol{u}_h^n , \boldsymbol{v}_h ) _\Sigma  = 
     \frac{\rho^s \epsilon }{\tau} ( \boldsymbol{ \dot{d}_h^{n-1} } + \tau \partial_\tau \boldsymbol{ \dot{d}_h^* } , \boldsymbol{v}_h )_\Sigma \\
    + \rho^f ( \partial_\tau \boldsymbol{u}_h^*, \boldsymbol{\mathcal{L}_h v_h} )_{\Omega^f} + a( \boldsymbol{u}_h^* , \boldsymbol{\mathcal{L}_h v_h} ) + b(p_h^*, \boldsymbol{ \mathcal{L}_h v_h } ) + l(\boldsymbol{v_h} )  \\
\end{aligned}
\right.
\end{eqnarray}

for all $ ( \boldsymbol{u}_h , q_h ) \in \boldsymbol{V}_h \times \boldsymbol{Q}_h \text{\quad with \quad } \boldsymbol{v}_h |_\Sigma \in \boldsymbol{W}_h $.

\vspace{1em}
2. Solid step: Find $ ( \boldsymbol{ \dot{d}_h^n } , \boldsymbol{d}_h^n ) \in \boldsymbol{W}_h \times \boldsymbol{W}_h $, such that 

\begin{eqnarray}    \label{eq:ERN5solid}
\left\{
\begin{aligned}
\boldsymbol{ \dot{d}_h^n }  =  \partial_\tau \boldsymbol{d}_h^n   \\
\rho^s \epsilon ( \partial_\tau  \boldsymbol{ \dot{d}^n } , \boldsymbol{w}_h )_\Sigma + \boldsymbol{a}^e ( \boldsymbol{d}_h^n , \boldsymbol{w}_h ) + \boldsymbol{a}^v ( \boldsymbol{ \dot{d} }_h^n , \boldsymbol{w}_h )  \\
    = -\rho^f ( \partial_\tau \boldsymbol{u}_h^n, \boldsymbol{\mathcal{L}_h w_h} ) - \boldsymbol{a} ( \boldsymbol{u}_h^n , \boldsymbol{\mathcal{L}_h w_h ) } -  \boldsymbol{b} ( p_h^n, \boldsymbol{\mathcal{L}_h w_h} ) 
\end{aligned}
\right.
\end{eqnarray}

for all $ \boldsymbol{w}_h \in \boldsymbol{W}_h $.
\par\noindent\rule{\textwidth}{.5pt}

\section{Application of jagged-time-step technique}
\subsection{New notations}
Let the expression $ \textbf{solveFluid} ( ( \boldsymbol{u}_h^n, p_h^n ) ; \quad (\boldsymbol{u}_h^{n-1}, p_h^{n-1} ), \tau_f, (\boldsymbol{ \dot{d}_h^\textit{m} }, \boldsymbol{d}_h^m ), extr = r )  $ denote the procedure solving the fluid part at  fluid time step $ n $ with known data of fluid from fluid time step $ n-1 $ and data of structure from step $ m $, where $ ( \boldsymbol{u}_h^n, p_h^n ) $ are the unknowns,  $ (\boldsymbol{u}_h^{n-1}, p_h^{n-1} ) $ are known fluid velocity and pressure from fluid time step $ n-1 $, $\tau_f $ is the length of fluid time step, $ (\boldsymbol{ \dot{d}_h^\textit{m} }, \boldsymbol{d}_h^m ) $ are the knownstructurevelocity and displacement fromstructuretime step $ m $, $ extr $ stands for the order of extrapolation forstructurevelocity and displacement ( see section $ 3.1 $ in \cite{rn} for details ). 

Analogously, $ \textbf{solveSolid} ( (\boldsymbol{ \dot{d}_h^\textit{m} }, \boldsymbol{d}_h^m ); \quad (\boldsymbol{ \dot{d}_h^\textit{m-1} }, \boldsymbol{d}_h^{m-1} ) , \tau_s,   ( \boldsymbol{u}_h^n, p_h^n ) , extr = r) $ denotes the procedure solving thestructurepart atstructuretime step $ m $ with known data of structure from step $ m-1 $ and data of fluid from fluid time step $ n $.

\subsection{Explicit Robin-Neumann scheme rewritten}
With the above notations, Fernandez's \textbf{ Explicit Robin-Neumann  scheme (intrinsic formulation)} can be rewritten as 

\par\noindent\rule{\textwidth}{1pt}
\textbf{(Fernandez) Explicit Robin-Neumann  scheme (intrinsic formulation)} 
\par\noindent\rule{\textwidth}{.5pt}

\noindent
Given final time $ T_{final} $, time step length $ \tau_3 $. 

\noindent
Let $  N = \frac{ T_{final}  }{\tau_3} $.

\noindent
\textbf{For} $ n = 1, 2, ..., N $, \textbf{do}

\indent
$ \textbf{solveFluid} ( ( \boldsymbol{u}_h^n, p_h^n ) ; \quad (\boldsymbol{u}_h^{n-1}, p_h^{n-1} ), \tau_3, (\boldsymbol{ \dot{d}_h^\textit{n-1} }, \boldsymbol{d}_h^{n-1} ), extr = r )  $ ;

\indent
$ \textbf{solveSolid} ( (\boldsymbol{ \dot{d}_h^\textit{n} }, \boldsymbol{d}_h^n ); \quad (\boldsymbol{ \dot{d}_h^\textit{n-1} }, \boldsymbol{d}_h^{n-1} ) , \tau_3,   ( \boldsymbol{u}_h^n, p_h^n ), extr = r )  $;
  
\noindent
\textbf{end for}
\par\noindent\rule{\textwidth}{1pt}

\subsection{Explicit Robin-Neumann scheme with jagged-time-step technique}

\par\noindent\rule{\textwidth}{1pt}
\textbf{Algorithm 1} Explicit Robin-Neumann scheme with jagged-time-step technique
\par\noindent\rule{\textwidth}{.5pt}

\noindent
Given final time $ T_{final} $, coarse time step length $ \tau_{coarse} = 10*\tau_{3}$, number of fluid steps within each coarse time interval $ N_f $, number of structure steps within each coarse time interval $ N_s $, $ N_f < 10 $, $ N_f + N_s \leq 20 $.

\noindent
Let $ N_{coarse} = \frac{ T_{final} }{ \tau_{coarse} }, \tau_f = \frac{ \tau_{coarse} }{ N_f }, \tau_s = \frac{ \tau_{coarse} }{ N_s } , n_{global} = 1 $.

\noindent
\textbf{For} $ i = 1, 2, ... , N_{coarse} $, within the coarse time interval $ ( (i-1)*\tau_{coarse}, i*\tau_{coarse} ] $, \textbf{do}

\indent
\textbf{For} $ ( \textbf{integer } m =1; m \leq N_s; m++ ) $, \textbf{do}

\indent \indent
\textbf{For} $ ( \textbf{integer } n = n_{global}; n \leq N_f; n++ ) $, \textbf{do}

\indent \indent \indent
\textbf{If} $ (m-1)*\tau_s < n*\tau_f \leq m*\tau_s $, \textbf{do}

\indent \indent \indent\indent  
$ \textbf{solveFluid} ( ( \boldsymbol{u}_h^n, p_h^n ) ; \quad (\boldsymbol{u}_h^{n-1}, p_h^{n-1} ), \tau_f, (\boldsymbol{ \dot{d}_h^\textit{m-1} }, \boldsymbol{d}_h^{m-1} ), extr = r )  $ ; \\
\indent \indent \indent\indent  
$ n_{global} = n + 1 $ ;     // remember past steps

\indent\indent\indent
\textbf{else do}

\indent \indent \indent\indent  
\textbf{break} the \textbf{for} loop;

\indent\indent\indent
\textbf{end if}

\indent\indent
\textbf{end for}

\indent\indent
$ \textbf{solveSolid} ( (\boldsymbol{ \dot{d}_h^\textit{m} }, \boldsymbol{d}_h^m ); \quad (\boldsymbol{ \dot{d}_h^\textit{m-1} }, \boldsymbol{d}_h^{m-1} ) , \tau_s,   ( \boldsymbol{u}_h^{n_{global} -1 }, p_h^{n_{global} - 1}  ) , extr = r) ; $

\indent
\textbf{end for}  \\
\indent
$ n_{global} = 1 $;

\noindent
\textbf{end for}
\par\noindent\rule{\textwidth}{.5pt}

\subsection{Some instances of Algorithm 1}
\begin{figure} [h]
\centering
\includegraphics[width = 0.99\textwidth]{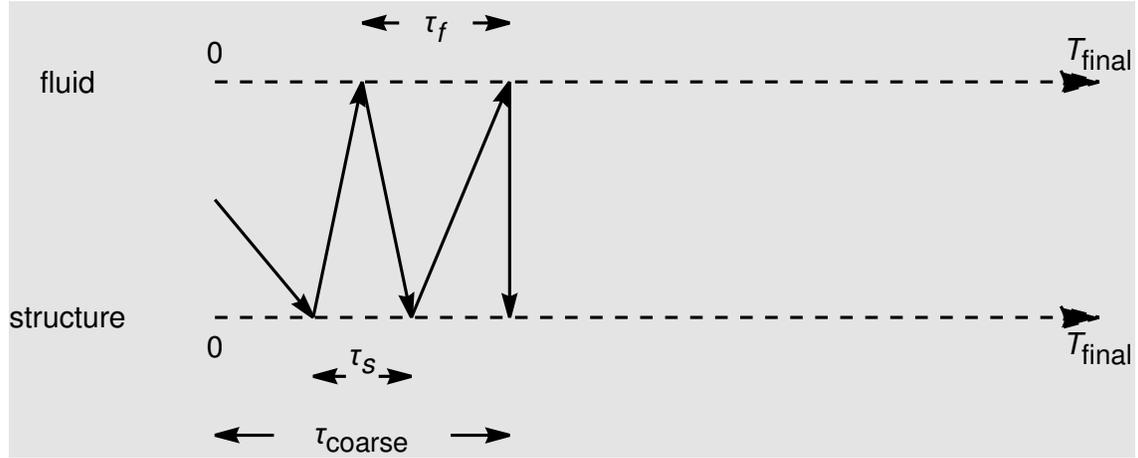}
\caption{"F 2 S 3" Algorithm 1}
\label{f2_s3}
\end{figure}

\begin{figure} [h]
\centering
\includegraphics[width = 0.99\textwidth]{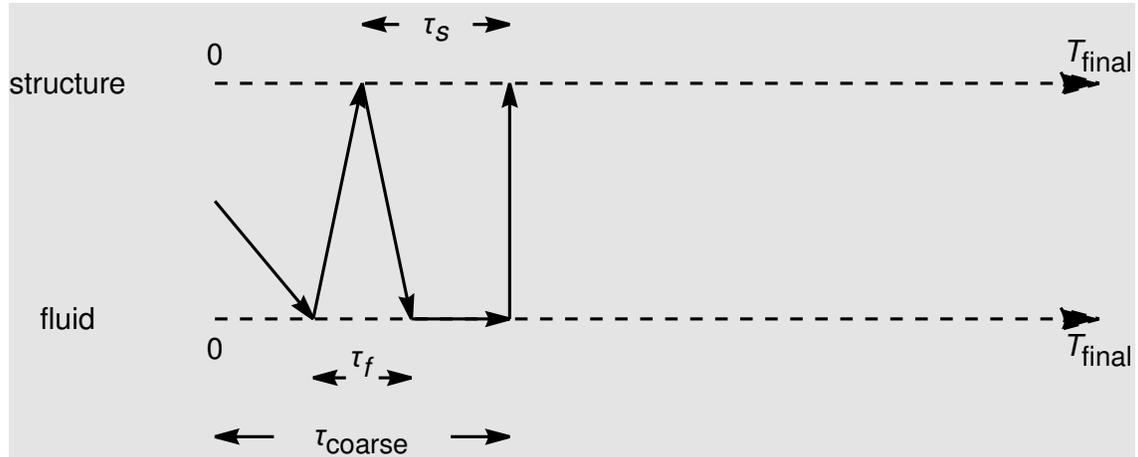}
\caption{"F 3 S 2" Algorithm 1}
\label{f3_s2}
\end{figure}

For short, let the expression "'F $ N_f $ S $ N_s $' Algorithm 1 " stand for \textbf{Algorithm 1} with given $ N_f $ and $ N_s $. Figure \ref{f2_s3} and \ref{f3_s2} describe the procedure of "'F 2 S 3' Algorithm 1" and "'F 3 S 2' Algorithm 1", respectively. Arrows therein point to the time when either the fluid or the structure is solved. 

\subsection{Two special cases of Algorithm 1}
The two cases of \textbf{Algorithm 1} with $ N_f = 1 $ or $ N_s = 1 $ are of special interests because they directly adopt the multiple-time-technique from \cite{multirate}. They are named \textbf{Algorithm 1.1} and \textbf{Algorithm 1.2}, respectively.

\par\noindent\rule{\textwidth}{1pt}
\textbf{Algorithm 1.1} A special case of \textbf{Algorithm 1}
\par\noindent\rule{\textwidth}{.5pt}
\noindent
Set $ N_f = 1 $ in \textbf{Algorithm 1}

\par\noindent\rule{\textwidth}{.5pt}
\vspace{2em}

\par\noindent\rule{\textwidth}{1pt}
\textbf{Algorithm 1.2} The other special case of \textbf{Algorithm 1} 
\par\noindent\rule{\textwidth}{.5pt}
\noindent
Set $ N_s = 1 $ in \textbf{Algorithm 1}
\par\noindent\rule{\textwidth}{.5pt}

\section{Numerical experiments}

\subsection{Configuration }
\begin{figure} [h]
\centering
\includegraphics[width = 0.99\textwidth]{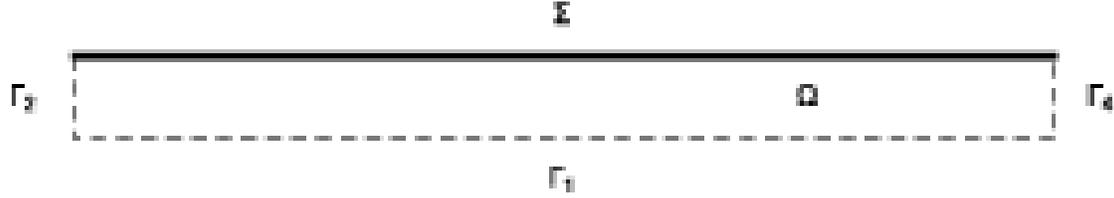}
\caption{Geometrial configuration}
\label{domain}
\end{figure}

The test-case as Section \textbf{6.1} in \cite{rn} is adopted,  except that in (\ref{eq:solid})  set  \\ 
$$
\boldsymbol{L^v \dot{d} } =
\begin{pmatrix}
0 \\
0
\end{pmatrix} . 
$$
Everything else remains intact. Namely, the fluid is defined on $ \Omega = [0, L] \times [0, R] $, where $ L = 6, R = 0.5 $ (all the quantities are under CGS system) , with $ \partial \Omega =  \Gamma_1 \cup \Gamma_2 \cup \Sigma \cup \Gamma_4 $ (see Figure \ref{domain}  )   . 

\noindent
On $ \Gamma_2 $ a sinusoidal pressure $ P (t) = P_{max} ( 1 - cos(2t \pi / T^\star )) /2  $ is prescribed, with $ P_{max} = 2*10^4 $  when $ 0 \leq t \leq T^\star $, $ P_{max} = 0 $ when $ t > T^\star , T^\star = 5*10^{-3} $. Zero pressure is imposed on $ \Gamma_4 $ and a symmetry condition is applied on $ \Gamma_1 $. Physical parameters for the fluid are 

$$ \rho^f = 1.0, \quad \mu = 0.035 .$$ 

The structure is assumed to be a generalized string defined on $ \Sigma $ with the two end points ( $ x = 0, L $ ) fixed, which therefore has lower dimension than the fluid. In (\ref{eq:solid}), take 

$$
\boldsymbol{d} = 
\begin{pmatrix} 
0  \\ 
\boldsymbol{d}_y  
\end{pmatrix},
\quad
\boldsymbol{L^e d} = 
\begin{pmatrix}
0   \\
-\lambda_1  \partial_{xx} \boldsymbol{d}_y  + \lambda_0 \boldsymbol{d}_y
\end{pmatrix},
$$

\noindent
with

$$
\lambda_1 \myeq \frac{E \epsilon}{2 (1+ \nu) } , 
\quad
\lambda_0 \myeq \frac{E \epsilon}{R^2 (1- \nu^2) } .
$$

\noindent
Physical parameters for thestructureare 

$$ E = 0.75 \times 10^6, \quad \epsilon = 0.1, \quad \nu = 0.5,  \quad \rho^s = 1.1 .  $$

All algorithms mentioned above are implemented using FreeFem++ \cite{ff++}. The Lagrange $ \textbf{P}_1 $ finite element is employed for both the fluid and structure, with symmetric pressure stabilization method introduced in \cite{stabilization} . The order of extrapolation is set to $ 1 $ ( $ extr = r = 1 $ ) . All run from time $ 0 $ to the final $ T_{final} = 0.015 (s) $.

To demonstrate both the $ h-$uniformity and the order of convergence in time, the time and space are refined at the same rate. Particularly, in Fernandez's \textbf{ Explicit Robin-Neumann  scheme (intrinsic formulation)}, set 

$$  (\tau_3, h ) = \frac{ (5*10^{-4}, 0.1) }{ 2^{rate} }, \quad rate = 0, 1, 2, 3, 4 . $$

\noindent
In \textbf{Algorithm 1, 1.1} and \textbf{1.2}, set 

$$ (\tau_{coarse}, h ) = ( 10*\tau_3, h ) = \frac{ (5*10^{-3}, 0.1) }{ 2^{rate} }, \quad rate = 0, 1, 2, 3, 4 , $$

\subsection{Numerical results} \label{numerical_results}
By comparing solutions of the above scheme to the reference solution generated by a fully implicit scheme with high space-time grid resolution ( $ \tau = 10^{-6}, h = 3.125 \times 10^{-3}  $ ), relative errors $ E_{rate} $ $ ( rate = 0, 1, 2, 3, 4 ) $ in elastic energy-norm (see \cite{rn} )  and time-convergence orders $ O_{rate} $  ( $ rate = 1, 2, 3, 4 $ )  defined as

$$ O_{rate} = \frac{\textrm{Log} \frac{ E_{rate} }{ E_{rate-1} } }{ \textrm{Log} \frac{1}{2} } , $$
for structure displacement corresponding to different refinement rates are computed.

Table \ref{table:ern} displays numerical results of Fernandez's \textbf{Explicit Robin-Neumann}  scheme. As predicted by the theoretical analysis in \cite{rn}, time-convergence order  of this scheme approaches $ 1 $ as both time and space get refined and is expected to reaches $ 1 $ as the refinement continues.

Table \ref{table:some_IMERN} reports some instances of \textbf{Algorithm 1} with $ N_f + N_s = 20, N_f < 10 $.

Table \ref{table:some_other_IMERN} exhibits some other instances of \textbf{Algorithm 1} with $ N_f + N_s < 20, N_f < 10 $.

Tables above demonstrate that as refinement \textit{rate} increases, instances of \textbf{Algorithm 1} reported in Table \ref{table:some_IMERN} and Table \ref{table:some_other_IMERN} obtain higher and higher convergence orders. Starting from certain refinement rates, their convergence orders and accuracy become higher than those of Fernandez's \textbf{Explicit Robin-Neumann} scheme.

There are some other observations. For example,  with $ N_f $ fixed, increasing $ N_s $ generally leads to lower relative errors and higher time-convergence orders. For example,  "F 6 S 13" Algorithm 1 works slightly better than "F6 S 12" Algorithm 1 in accuracy and convergence orders.

Neither \textbf{Algorithm 1.1} nore \textbf{Algorithm 1.2} yield satisfactory results. Some instances are not even stable. Among the stable instances, "F 1 S 20" performs best, but still much worse than Fernandez's \textbf{Explicit Robin-Neumann}.  Table \ref{table:an_mern} reports results of "'F 1 S 20' Algorithm 1.1", namely "'F 1 S 20' Algorithm 1".

\begin{table}[h!]
\centering
\begin{tabular}{ l l l }
\hline
rate & $ E_{rate} $  & $ O_{rate} $    \\
\hline
0 &			 0.959089  &  \\
1 &			 0.719217 &	 0.415238 \\
2 &			 0.435036 &	 0.725292 \\
3 &			 0.241714 &	 0.847834 \\
4 &			 0.128601 &	 0.910399 \\
\hline
\end{tabular}
\caption{Numerical results of \textbf{Explicit Robin-Neumann}  scheme ("F 10 S 10" Algorithm 1) }
\label{table:ern}
\end{table}

\begin{table}[h!]
\centering
\begin{tabular}{ c c }

\begin{tabular}{ l l l }
"F 4 S 16" &  Algorithm 1  & \\
\hline
rate & $ E_{rate} $  & $ O_{rate} $    \\
\hline
0 &			 1.07566  &  \\
1 &			 1.03863 &	 0.0505403 \\
2 &			 0.673546 &	 0.624833 \\
3 &			 0.289224 &	 1.21959 \\
4 &			 0.112203 &	 1.36608 \\
\hline
\end{tabular}

&

\begin{tabular}{ l l l }
"F 5 S 15" & Algorithm 1 & \\
\hline
rate & $ E_{rate} $  & $ O_{rate} $    \\
\hline
0 &			 1.0712  &  \\
1 &			 0.949946 &	 0.17331 \\
2 &			 0.555675 &	 0.773604 \\
3 &			 0.24079 &	 1.20647 \\
4 &			 0.102342 &	 1.23438 \\
\hline
\end{tabular}

\\
& \\
& \\

\begin{tabular}{ l l l }
"F 6 S 14" & Algorithm 1 & \\
\hline
rate & $ E_{rate} $  & $ O_{rate} $    \\
\hline
0 &			 1.05564  &  \\
1 &			 0.882138 &	 0.259042 \\
2 &			 0.501674 &	 0.814254 \\
3 &			 0.231413 &	 1.11628 \\
4 &			 0.106829 &	 1.11517 \\
\hline
\end{tabular}

&

\begin{tabular}{ l l l }
"F 7 S 13" & Algorithm 1 & \\
\hline
rate & $ E_{rate} $  & $ O_{rate} $    \\
\hline
0 &			 1.03319  &  \\
1 &			 0.824216 &	 0.326011 \\
2 &			 0.465121 &	 0.825416 \\
3 &			 0.226727 &	 1.03665 \\
4 &			 0.110134 &	 1.0417 \\
\hline
\end{tabular}

\end{tabular}
\caption{Numerical results of some instances of Algorithm 1 ( $ N_f + N_s = 20, N_f < 10 $ ) }
\label{table:some_IMERN}
\end{table}

%
%
%
%

\begin{table}[h!]
\centering
\begin{tabular}{ c c }

\begin{tabular}{ l l l }
"F 4 S 14" & Algorithm 1 & \\
\hline
rate & $ E_{rate} $  & $ O_{rate} $    \\
\hline
0 &                      1.07127  &  \\
1 &                      1.04302 &       0.0385553 \\
2 &                      0.688734 &      0.598748 \\
3 &                      0.30437 &       1.17812 \\
4 &                      0.121974 &      1.31925 \\
\hline
\end{tabular}

&

\begin{tabular}{ l l l }
"F 4 S 15" & Algorithm 1 & \\
\hline
rate & $ E_{rate} $  & $ O_{rate} $    \\
\hline
0 &                      1.0744  &  \\
1 &                      1.0464 &        0.0380968 \\
2 &                      0.68863 &       0.603634 \\
3 &                      0.302521 &      1.18669 \\
4 &                      0.120325 &      1.3301 \\
\hline
\end{tabular}

\\

& \\
& \\

\begin{tabular}{ l l l }
"F 5 S 12" & Algorithm 1 & \\
\hline
rate & $ E_{rate} $  & $ O_{rate} $    \\
\hline
0 &                      1.07029  &  \\
1 &                      0.971898 &      0.139125 \\
2 &                      0.593959 &      0.710442 \\
3 &                      0.275873 &      1.10636 \\
4 &                      0.125772 &      1.13319 \\
\hline
\end{tabular}

& 

\begin{tabular}{ l l l }
"F 5 S 13" & Algorithm 1 & \\
\hline
rate & $ E_{rate} $  & $ O_{rate} $    \\
\hline
0 &			 1.07208  &  \\
1 &			 0.967564 &	 0.147984 \\
2 &			 0.584823 &	 0.726357 \\
3 &			 0.26601 &	 1.13652 \\
4 &			 0.118319 &	 1.1688 \\
\hline
\end{tabular}

\\
& \\ &  \\

\begin{tabular}{ l l l }
"F 5 S 14" & Algorithm 1 & \\
\hline
rate & $ E_{rate} $  & $ O_{rate} $    \\
\hline
0 &			 1.0723  &  \\
1 &			 0.964902 &	 0.152254 \\
2 &			 0.575995 &	 0.744326 \\
3 &			 0.257087 &	 1.1638 \\
4 &			 0.11213 &	 1.19708 \\
\hline
\end{tabular}

&

\begin{tabular}{ l l l }
"F 6 S 12" & Algorithm 1 & \\
\hline
rate & $ E_{rate} $  & $ O_{rate} $    \\
\hline
0 &			 1.05435  &  \\
1 &			 0.879988 &	 0.260798 \\
2 &			 0.504881 &	 0.80154 \\
3 &			 0.239873 &	 1.07367 \\
4 &			 0.114482 &	 1.06715 \\
\hline
\end{tabular}

\\ & \\ & \\

\begin{tabular}{ l l l }
"F 6 S 13" & Algorithm 1 & \\
\hline
rate & $ E_{rate} $  & $ O_{rate} $    \\
\hline
0 &			 1.05739  &  \\
1 &			 0.886937 &	 0.253604 \\
2 &			 0.512224 &	 0.792057 \\
3 &			 0.241241 &	 1.0863 \\
4 &			 0.11318 &	 1.09186 \\
\hline
\end{tabular}

&

\begin{tabular}{ l l l }
"F 7 S 12" & Algorithm 1 & \\
\hline
rate & $ E_{rate} $  & $ O_{rate} $    \\
\hline
0 &			 1.03817  &  \\
1 &			 0.834026 &	 0.315878 \\
2 &			 0.482387 &	 0.789901 \\
3 &			 0.242281 &	 0.99351 \\
4 &			 0.120704 &	 1.00521 \\
\hline
\end{tabular}

\end{tabular}
\caption{Numerical results of some other Algorithm 1 instances ( $ N_f + N_s < 20, N_f < 10 $ ) }
\label{table:some_other_IMERN}
\end{table}

\begin{table}[h!]
\centering
\begin{tabular}{ l l l }
\hline
rate & $ E_{rate} $  & $ O_{rate} $    \\
\hline
0 &			 1  &  \\
1 &			 1.04684 &	 -0.066041 \\
2 &			 1.12592 &	 -0.105063 \\
3 &			 1.07841 &	 0.0621985 \\
4 &			 0.639951 &	 0.752872 \\
\hline
\end{tabular}
  \caption{Numerical results of "F 1 S 20" Algorithm 1.1 (namely "F 1 S 20" Algorithm 1)}
\label{table:an_mern}
\end{table}

\subsection{Graphs of structure displacements}
\begin{figure}[h]
\centering
\subfloat{ \includegraphics[width = 0.5\textwidth]{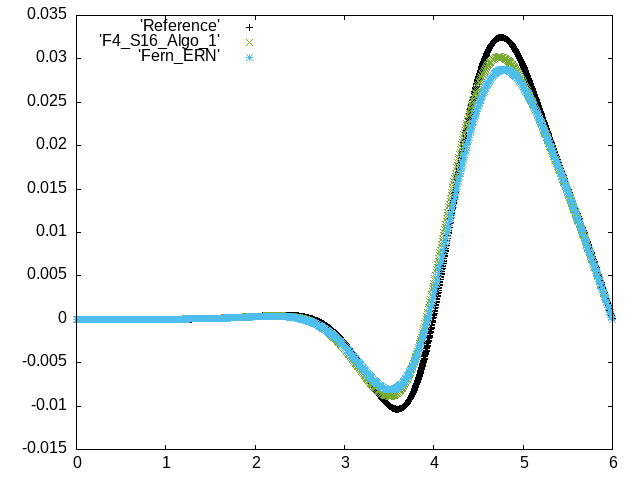} }
\subfloat{ \includegraphics[width = 0.5\textwidth]{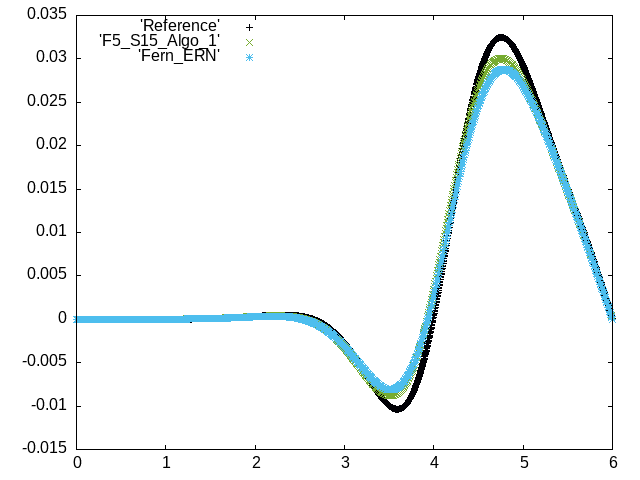} }   \\
\subfloat{ \includegraphics[width = 0.5\textwidth]{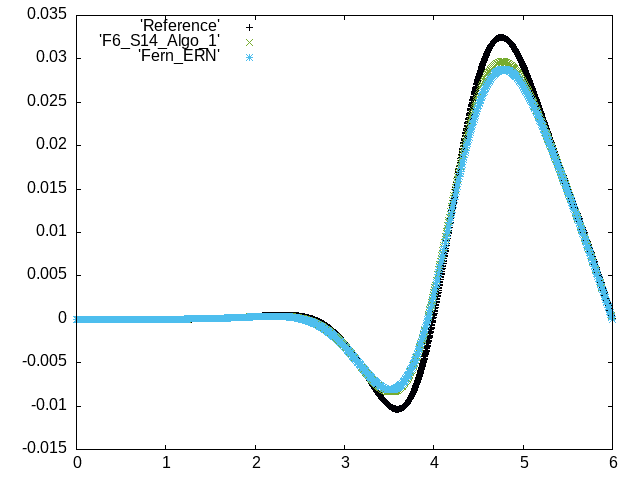} }
\subfloat{ \includegraphics[width = 0.5\textwidth]{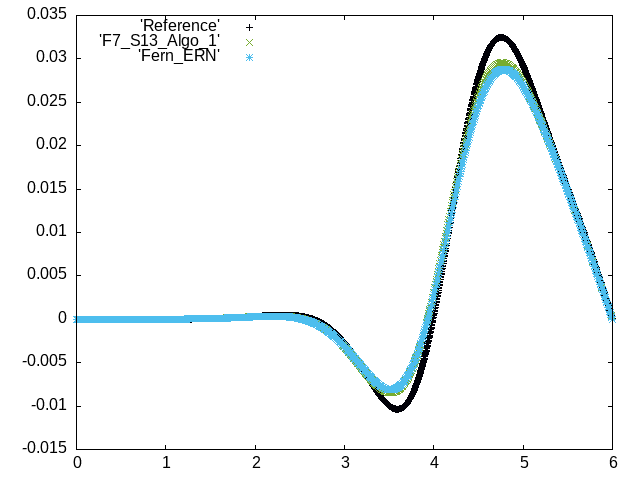} }
  \caption{Structure displacements of each instance in Table \ref{table:some_IMERN} , the reference and Fernandez's \textbf{Explicit Robin-Neumann} at $ T_{final} $}
\label{fig:some_IMERN}
\end{figure}

\begin{figure}[h]
\centering
\subfloat{ \includegraphics[width = 0.5\textwidth]{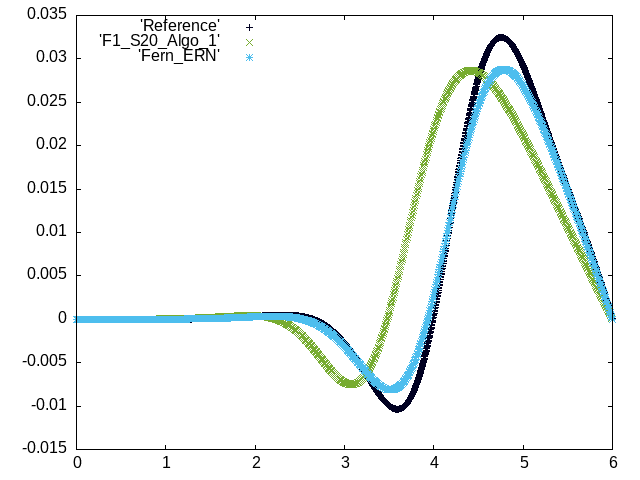} }
  \caption{ Structure displacements of the instance in Table \ref{table:an_mern} , the reference and Fernandez's \textbf{Explicit Robin-Neumann} at $ T_{final} $}
\label{fig:MERN}
\end{figure}

Figures \ref{fig:some_IMERN}, \ref{fig:some_other_IMERN} and \ref{fig:MERN} display structure displacements at time $ T_{final} $ and refinement $ rate = 4 $ for all instances presented in Tables \ref{table:some_IMERN}, \ref{table:some_other_IMERN} and \ref{table:an_mern} in comparison with those of the reference and Fernandez's \textbf{Explicit Robin-Neumann}. 


\begin{figure}[htb]
\centering
\subfloat{ \includegraphics[width = 0.5\textwidth]{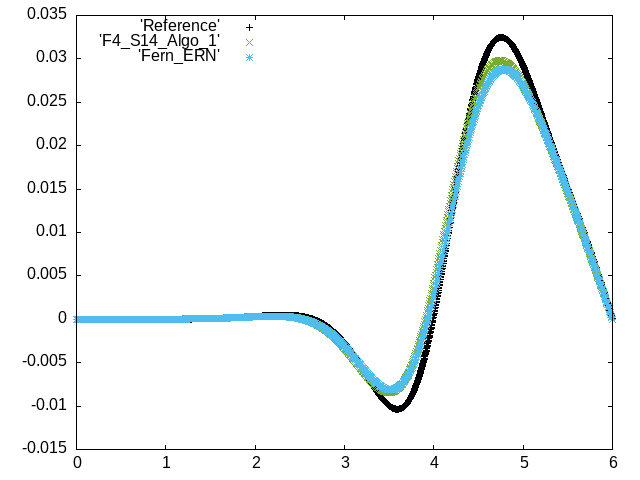} }
\subfloat{ \includegraphics[width = 0.5\textwidth]{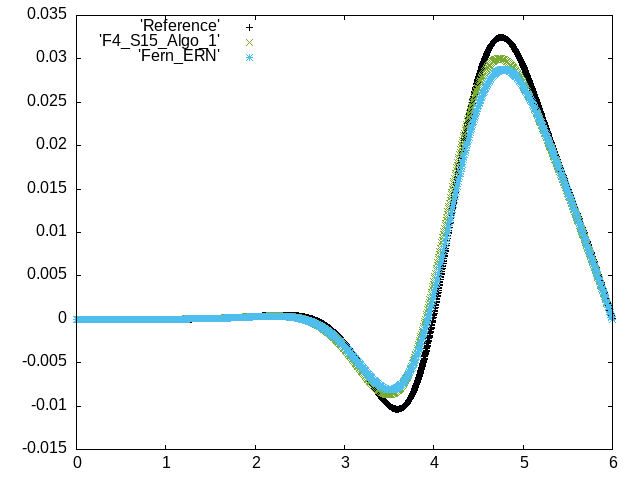} }   \\
\subfloat{ \includegraphics[width = 0.5\textwidth]{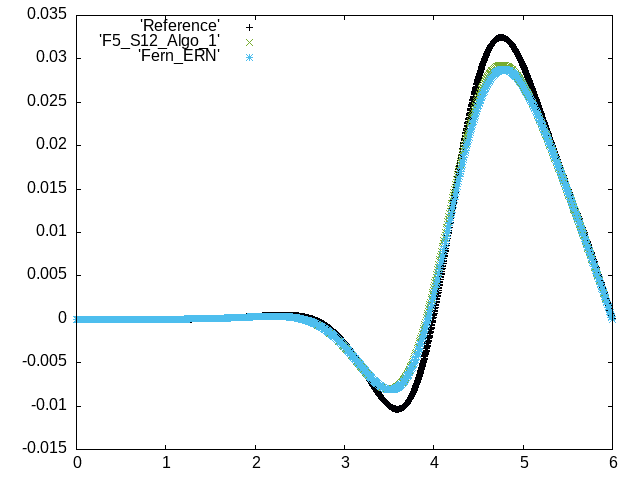} }
\subfloat{ \includegraphics[width = 0.5\textwidth]{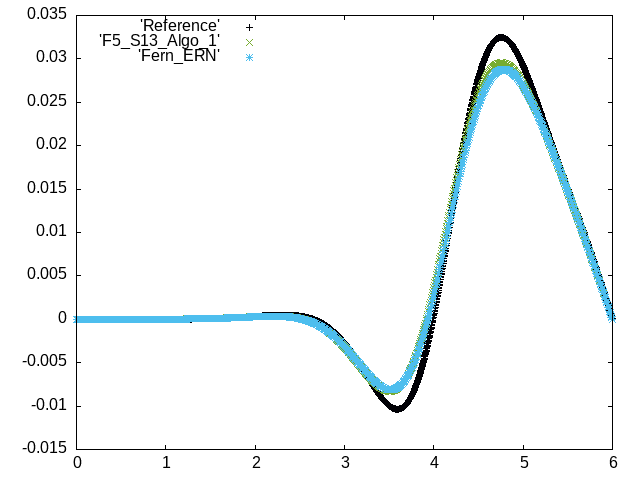} }   \\
\subfloat{ \includegraphics[width = 0.5\textwidth]{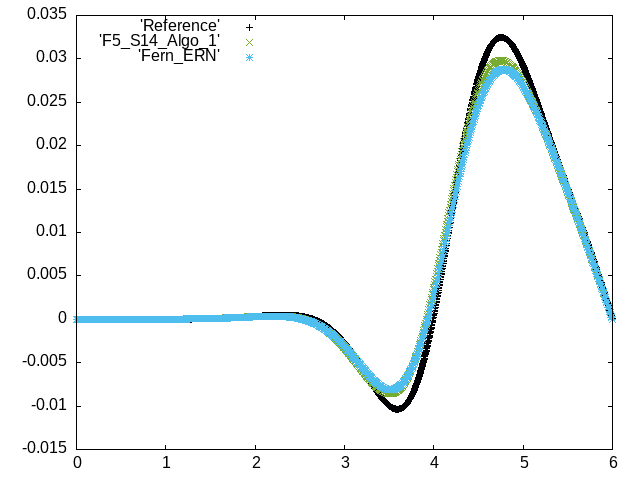} }
\subfloat{ \includegraphics[width = 0.5\textwidth]{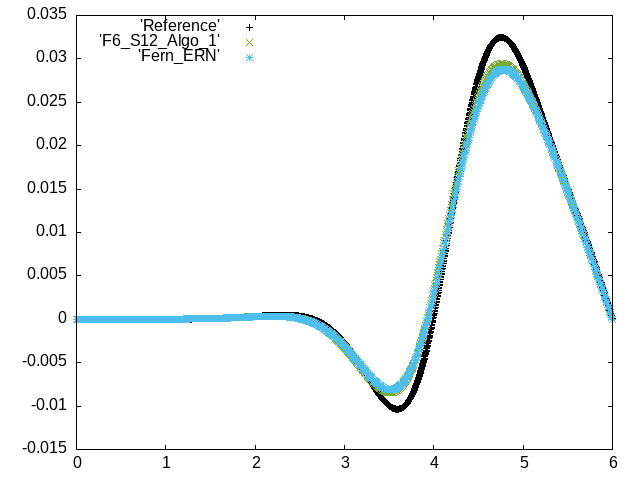} }   \\
  \caption{ Structure displacements of each instance in Table \ref{table:some_other_IMERN} , the reference and Fernandez's \textbf{Explicit Robin-Neumann} at $ T_{final} $}
\label{fig:some_other_IMERN}
\end{figure}

\begin{figure}[htb]  \ContinuedFloat  
\centering
\subfloat{ \includegraphics[width = 0.5\textwidth]{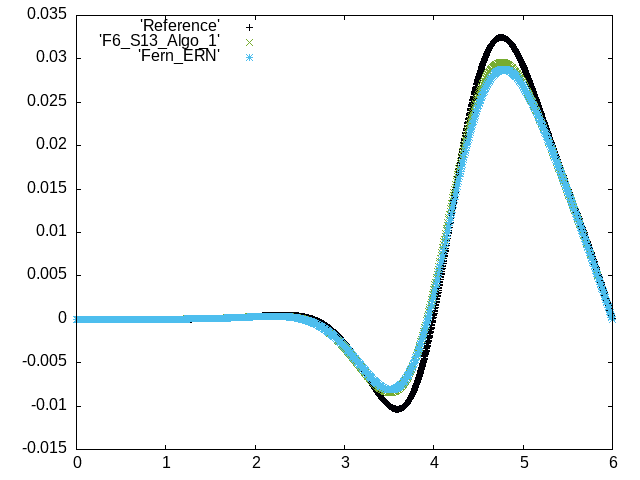} }
\subfloat{ \includegraphics[width = 0.5\textwidth]{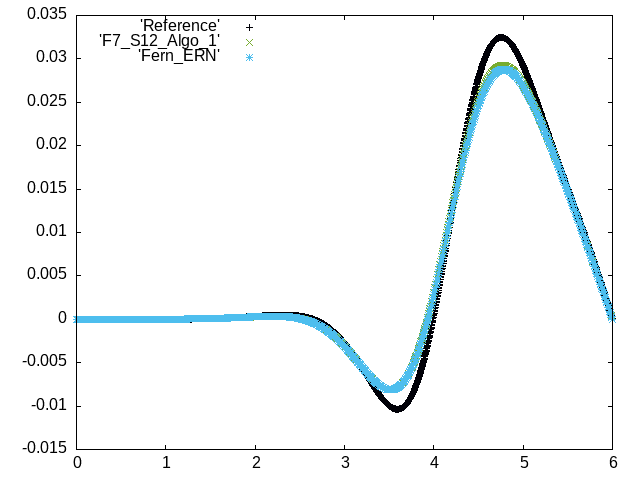} }   \\
  \caption{ Structure displacements of each instance in Table \ref{table:some_other_IMERN} , the reference and Fernandez's \textbf{Explicit Robin-Neumann} at $ T_{final} $}
\label{fig:some_other_IMERN}
\end{figure}

\section{Conclusions}
As refinement \textit{rate} increases, some instances of \textbf{Algorithm 1} obtain higher convergence orders and accuracy than the original \textbf{Explicit Robin-Neumann scheme} with lower cost. 

\section{Discussion and future work }
The ideas of jagged-time-step technique might applicable to other algorithms or problems.

\section{Acknowledgments}
This research did not receive any specific grant from funding agencies in the public, commercial, or not-for-profit sectors.

The ideas as well as algorithms described in this work (except those cited from other sources explicit stated), numerical experiments, results and conclusions are completed entirely independently, without any assistance from anyone else.

More than two years ago, Dr. Mingchao CAI found the article \cite{multirate} and suggested applying directly the multiple-time-technique therein to the coupling of incompressible fluid with thin-walled structure. Dr. Lian ZHANG adopted his ideas. However, this work has nothing to do with their ideas or work.

\clearpage
\newpage
\bibliographystyle{plain}
\bibliography{ref}  
\end{document}